\documentclass[10pt]{article}
\font\smallit=cmti10
\font\smalltt=cmtt10

\usepackage{amssymb,latexsym,amsmath,epsfig,amsthm} 
\usepackage[pdftex,pdfpagelabels]{hyperref}
\usepackage{mathtools}

\makeatletter

\renewcommand\section{\@startsection {section}{1}{\z@}
{-30pt \@plus -1ex \@minus -.2ex}
{2.3ex \@plus.2ex}
{\normalfont\normalsize\bfseries\boldmath}}

\renewcommand\subsection{\@startsection{subsection}{2}{\z@}
{-3.25ex\@plus -1ex \@minus -.2ex}
{1.5ex \@plus .2ex}
{\normalfont\normalsize\bfseries\boldmath}}

\renewcommand{\@seccntformat}[1]{\csname the#1\endcsname. }

\makeatother

\newtheorem{theorem}{Theorem}
\newtheorem{lemma}{Lemma}
\newtheorem{proposition}{Proposition}

\newtheorem{problem}{Problem}

\theoremstyle{definition}

\newtheorem{remark}{Remark}

\renewcommand\P{\mathbb{P}}
\newcommand\E{\mathbb{E}}

\newcommand\N{\mathbb{N}}

\newcommand\Log{{\operatorname{Log}}}
\newcommand\lcm{{\operatorname{lcm}}}
\renewcommand\gcd{{\operatorname{gcd}}}
\newcommand\eps{\varepsilon}

\begin{document}

\begin{center}
  \uppercase{\bf \boldmath Dense sets of natural numbers with unusually large least common multiples}
  \vskip 20pt
  {\bf Terence Tao}\\
  {\smallit Department of Mathematics, UCLA, Los Angeles, CA}\\
  {\tt tao@math.ucla.edu}\\
  \end{center}
  \vskip 20pt
  \centerline{\smallit Received: , Revised: , Accepted: , Published: } 
  \vskip 30pt

  \centerline{\bf Abstract}
  \noindent
  For any constant $C_0>0$, we construct a set $A \subset \N$ such that one has
  $$ \sum_{n \in A: n \leq x} \frac{1}{n} = \exp\left(\left(\frac{C_0}{2}+o(1)\right) (\log\log x)^{1/2} \log\log\log x \right)$$
  and
  $$ \sum_{n,m \in A: n, m \leq x} \frac{1}{\lcm(n,m)} \ll (e^{C_0^2}-1+o(1))  \left(\sum_{n \in A: n \leq x} \frac{1}{n}\right)^2$$
  as $x \to \infty$, with the growth rate given here optimal up to the dependence on $C_0$.   This answers in the negative a question of Erd\H{o}s and Graham, and also clarifies the nature of certain ``mostly coprime'' sets studied by Bergelson and Richter.

  \pagestyle{myheadings}
  \markright{\smalltt INTEGERS: 24 (2024)\hfill}
  \thispagestyle{empty}
  \baselineskip=12.875pt
  \vskip 30pt


\section{Introduction}

In this paper we resolve (in the negative, and in a near-optimal fashion) a question of Erd\H{o}s and Graham \cite{erdos} (problem \#442 in \url{https://erdosproblems.com}). Our analysis also clarifies the nature of certain ``mostly coprime'' sets studied by Bergelson and Richter \cite{bergelson}.

We begin with the former question. It is convenient to introduce the iterated logarithm notation\footnote{In this paper the iterated logarithm takes precedence over all other operations except for parentheses, thus for instance $\Log^{1/2}_2 x \Log_3 x$ is equal to $(\Log_2 x)^{1/2} (\Log_3 x)$ rather than $(\Log_2 (x \Log_3 x))^{1/2}$).}
$$ \Log x \coloneqq \max(\log x, 1); \quad \Log_2 x \coloneqq \Log\Log x; \quad \Log_3 x \coloneqq \Log\Log\Log x$$
for $x > 0$.

\begin{problem}[Erd\H{o}s \#442]\label{442}\cite[p. 88]{erdos}  Let $A \subset \N$ be such that
\begin{equation}\label{agrow}
   \frac{1}{\Log_2 x} \sum_{n \in A: n \leq x} \frac{1}{n} \to \infty
\end{equation}
as $x \to \infty$.  Is it true that
\begin{equation}\label{conc}
 \frac{1}{(\sum_{n \in A: n \leq x} \frac{1}{n})^2} \sum_{n,m \in A: n < m \leq x} \frac{1}{\lcm(n,m)} \to \infty
\end{equation}
as $x \to \infty$?
\end{problem}

Informally, Problem \ref{442} is asking how large a set $A$ can become assuming that the least common multiples $\lcm(n,m)$ are unusually large in a certain (exotic) averaged sense.  In view of the well-known identity
\begin{equation}\label{gcd}
  \lcm(n,m) = \frac{nm}{\gcd(n,m)},
\end{equation}
one can also interpret Problem \ref{442} as studying the largest size of sets $A$ for which the greatest common divisors $\gcd(n,m)$ of elements $n,m \in A$ are unusually small on average.  See also \cite{dmt} for some similar problems of this type regarding minimization of G\'al-type sums such as $\sum_{n,m \in A: n,m \leq x} \frac{\gcd(n,m)}{m+n}$.

\begin{remark} At the opposite extreme, sets with unusually \emph{large} greatest common divisors (or unusually small least common multiples) have recently been studied in \cite{km}, \cite{gw}, as they are relevant to the resolution of the Duffin--Schaeffer conjecture in \cite{km}.  The aforementioned reference \cite[p. 87]{erdos} also states a problem in this context (problem \#441 in \url{https://erdosproblems.com}, see also problems B26, E2 in \cite{guy}), regarding the size of the largest set $A$ of integers for which $\lcm(n,m) \leq x$ for all $n,m \in A$, which was subsequently solved rather satisfactorily in \cite{chen, chen-2, chen-3}; the answer is $\sqrt{9x/8}-2$, plus a non-negative error of size $O(x^{1/2} \Log^{-1/2} x \Log_2 x)$ that goes to infinity (by at least the inverse tower exponential function) for an infinite sequence of $x$.
\end{remark}

The presence of the $\Log_2 x$ growth rate in Problem \ref{442} was likely motivated by the example of $A$ being the set of prime numbers, since a straightforward application of Mertens' theorem shows that both expressions appearing in this problem are comparable to $1$ for large $x$ in this case (note that $\lcm(p_1,p_2)=p_1 p_2$ for any distinct primes $p_1,p_2$).  The potential relevance of this example is also supported by the easy observation that the maximal value of $\sum_{n \in A: n \leq x} \frac{1}{n}$ if the elements of $A$ are assumed to be pairwise coprime (which is an evident way to prevent \eqref{conc} from occurring) is attained when $A$ is the set of prime numbers (together with $1$).

Observe that the diagonal expression
$$ \frac{1}{(\sum_{n \in A: n \leq x} \frac{1}{n})^2} \sum_{n,m \in A: n = m \leq x} \frac{1}{\lcm(n,m)}$$
goes to zero under the hypotheses of this theorem, since $\lcm(n,n)=n$.  By symmetry, the conclusion \eqref{conc} of this problem is then equivalent to
$$ \frac{1}{(\sum_{n \in A: n \leq x} \frac{1}{n})^2} \sum_{n,m \in A: n, m \leq x} \frac{1}{\lcm(n,m)} \to \infty.$$
Using \eqref{gcd}, we can write this as
$$ \frac{\sum_{n,m \in A: n, m \leq x} \frac{\gcd(n,m)}{nm}}{\sum_{n,m \in A: n, m \leq x} \frac{1}{nm}} \to \infty.$$
One can interpret this claim probabilistically: if we let $\mathbf{n} = \mathbf{n}_x$ denote a random element of $\{n \in A: n \leq x\}$ drawn with the logarithmic probability distribution
\begin{equation}\label{n-prob}
  \mathbb{P}( \mathbf{n} = n ) = \frac{1/n}{\sum_{m \in A: m \leq x} \frac{1}{m}},
\end{equation}
and let $\mathbf{m} = \mathbf{m}_x$ denote an independent copy of $\mathbf{n}$, then the claim \eqref{conc} is now equivalent to the assertion that
\begin{equation}\label{gcd-grow}
   \E \gcd(\mathbf{n}, \mathbf{m}) \to \infty.
\end{equation}
Note that as $\gcd(\mathbf{n}, \mathbf{m}) \geq 1$, we trivially have
\begin{equation}\label{trivial}
  \E \gcd(\mathbf{n}, \mathbf{m}) \geq 1
\end{equation}
and so (by reversing the above manipulations) the left-hand side of \eqref{conc} must be $\gg 1$ for sufficiently large $x$. Informally, Problem \ref{442} asks whether any set that is significantly denser than the primes (in a logarithmic sense) must start having large pairwise greatest common divisors on average.  A key advantage of the reformulation \eqref{gcd-grow} of the problem is that one can now use the standard Gauss identity
\begin{equation}\label{gcd-standard}
  \gcd(n,m) = \sum_{d|n,m} \phi(d),
\end{equation}
where $\phi$ is the Euler totient function, together with the fact that $\mathbf{m}$ is an independent copy of $\mathbf{n}$, to write the left-hand side of \eqref{gcd-grow} as a weighted sum of squares:
\begin{equation}\label{sum-square}
  \E \gcd(\mathbf{n}, \mathbf{m}) = \sum_d \phi(d) \P(d|\mathbf{n})^2.
\end{equation}
Thus, the claim \eqref{conc} (or \eqref{gcd-grow}) can also be interpreted as an assertion that typical elements of $A$ do not have unusually few divisors. Note that the trivial lower bound \eqref{trivial} corresponds to the $d=1$ term in the right-hand side of \eqref{sum-square}.

The defect\footnote{We thank Joel Moreira for this observation.}
\begin{equation}\label{defect}
  \E \gcd(\mathbf{n}, \mathbf{m}) - 1
\end{equation}
in the inequality \eqref{trivial} was studied in \cite{bergelson} , where it was denoted as
$$\E_{m \in A; m \leq x}^{\log} \E_{n \in A; n \leq x}^{\log} \Phi(n,m),$$
and can be viewed as a (weighted) measure of how often two elements of $A$ of size at most $x$ share a common factor.
By \eqref{sum-square}, this quantity can also be written as
\begin{equation}\label{defect-alt}
  \sum_{d>1} \phi(d) \P(d|\mathbf{n})^2.
\end{equation}
The authors of \cite{bergelson} were particularly interested in obtaining sets $A$ for which the defect \eqref{defect} was small (less than $\eps$ for a given fixed small $\eps>0$), as this allowed for a useful inequality of Tur\'an--Kubilius type to be available for such sets $A$, which had applications to generalizations of the Sarnak conjecture \cite{sarnak}; see in particular \cite[Proposition 2.1]{bergelson} and the surrounding discussion.

In the case where $A$ consists of the primes, one can verify from Mertens' theorem that $\P(d|\mathbf{n})$ is asymptotic to $\frac{1}{d \Log_2 d}$ when $d$ is prime (or equal to $1$), and vanishes otherwise.  Hence the sum \eqref{sum-square} converges with a ``doubly logarithmic amount'' of room to spare; indeed, from another application of Mertens' theorem we see that the defect \eqref{defect} decays like\footnote{In this paper we use the usual asymptotic notation $X \ll Y$, $Y \gg X$ or $X = O(Y)$ to denote $|X| \leq CY$ for an absolute constant $C$, $X \asymp Y$ to denote $X \ll Y \ll X$, and $o(1)$ to denote a quantity that goes to zero as $x \to \infty$.}  $(1+o(1))/\Log_2 x$.  This suggests that the primes are in fact not the densest set for which \eqref{conc} fails.  For instance, if one instead takes $A$ to be the set of squarefree numbers with exactly $k$ prime factors for any fixed $k \geq 1$, a standard calculation again based on Mertens' theorem (see also Lemma \ref{lgr} below) shows that the left-hand side of \eqref{agrow} now grows like $(1+o(1)) \Log^{k-1}_2 x / k!$, but the left-hand side of \eqref{conc} stays bounded in the limit $x \to \infty$ (and the defect \eqref{defect} decays like $(1+o(1)) k^2/\Log_2 x$).  In particular, the answer to Problem \ref{442} is negative; this result and construction was already implicitly observed in \cite{bergelson} (see the discussion after \cite[Proposition 2.1]{bergelson}), although the authors seem to have been unaware of Problem \ref{442}. In the remarks after \cite[(2.5)]{bergelson}, it is also noted that the situation appears to be significantly worse if one performs natural averaging rather than logarithmic averaging for these problems, and we do not pursue that variant of the problem here.

In this paper, we can identify (up to constants) the optimal (logarithmic) growth rate of $A$ for which \eqref{conc} stays bounded.  We first present a simplified version of the main theorem:

\begin{theorem}[Negative answer to \#442]\label{simplified-main} For any $C_0 > 0$, there exists a set $A$ of natural numbers such that
\begin{equation}\label{442-1}
 \sum_{n \in A: n \leq x} \frac{1}{n} = \exp\left( \left(\left(\frac{C_0}{2}+o(1)\right) \Log_2^{1/2} x \Log_3 x \right) \right)
\end{equation}
and
\begin{equation}\label{442-2}
   \sum_{n, m \in A: n, m \leq x} \frac{1}{\lcm(n,m)} \ll_{C_0} \left( \sum_{n \in A: n \leq x} \frac{1}{n} \right)^2
\end{equation}
as $x \to \infty$.  Up to the choice of constant $C_0$, the growth rate in \eqref{442-1} is otherwise optimal for sets that obey \eqref{442-2}.
\end{theorem}

In other words, the optimal rate of growth of $\sum_{n \in A: n \leq x} \frac{1}{n}$ in order to avoid the conclusion \eqref{conc} of Problem \ref{442} is not $\asymp \Log_2 x \asymp \exp( \Log_3 x)$ as stated in that problem, but is instead $\exp( O(\Log^{1/2}_2 x \Log_3 x ))$.

In fact, we prove a more precise statement. Let $C_0>0$ be a fixed parameter; this quantity could be either larger, smaller, or equal to $1$, but we recommend that the reader focus on the $C_0=1$ case on a first reading.  For any natural number $k \geq C_0$, let $x_k = x_{k,C_0}$ denote the real number $x_k \coloneqq \exp\exp(k^2/C_0^2)$ (so that $k = C_0 \Log_2^{1/2} x_k$), and for any $x_k < x \leq x_{k+1}$, so that
$$\frac{k^2}{C_0^2} < \Log_2 x \leq \frac{(k+1)^2}{C_0^2} = \frac{k^2}{C_0^2} + O\left(\frac{k}{C_0^2}\right),$$
let $\psi(x) = \psi_{C_0}(x)$ denote the quantity
\begin{equation}\label{psi-def}
  \psi(x) \coloneqq 1 + \frac{h(x)^2}{h(x_{k+1})}
\end{equation}
with $h(x) = h_{C_0}(x)$ defined by
\begin{equation}\label{h-def}
h(x) \coloneqq \Log_2 x - \frac{k^2}{C_0^2} + 1
\end{equation}
(so in particular $h(x_{k+1}) \asymp k/C_0^2$ for large $k$), and let
$F(x) = F_{C_0}(x)$ denote the quantity
\begin{equation}\label{F-def}
   F(x) \coloneqq \psi(x) \frac{k^{2k} C_0^{-2k}}{k!},
\end{equation}
thus $F$ and $\psi$ are piecewise quadratic left-continuous functions of $\Log_2 x$, defined for $x$ large enough, and $h$ is similarly piecewise linear in $\Log_2 x$.  Observe that on the interval $x_k < x \leq x_{k+1}$ for large $k$, $\psi$ and $h$ both increase from $\lim_{x \to x_k^+} \psi(x) = \lim_{x \to x_k^+} h(x) = 1$ to $\psi(x_{k+1}) \asymp h(x_{k+1}) \asymp k/C_0^2$, and as a consequence $F$ is ``continuous up to constants'' in the sense that
\begin{equation}\label{F-cts}
  F(x_k) \asymp \frac{k}{C_0^2} \frac{(k-1)^{2(k-1)} C_0^{-2(k-1)}}{(k-1)!} \asymp \frac{k^{2k} C_0^{-2k}}{k!} = \lim_{x \to x_k^+} F(x)
\end{equation}
for sufficiently large $k$.  For technical reasons (that mainly concern the case when $C_0$ is small) it will be important to observe that $\psi$ is ``flat'' for some distance to the right of $x_k$, and specifically that
\begin{equation}\label{psi-flat}
  \psi(x) \asymp 1
\end{equation}
whenever $x_k < x \leq x_{k+1}$ with $h(x) \ll k^{1/2}/C_0$ and $k$ large; in particular we have
\begin{equation}\label{F-flat}
  F(x) \asymp F(x_k)
\end{equation}
in the same regime.  This is the main reason why we select a quadratic form \eqref{psi-def} as the first factor in the right-hand side of \eqref{F-def}, rather than the more natural linear form $h(x)$.  From Stirling's formula (and the fact that $k \asymp C_0 \Log_2^{1/2} x$) one can easily establish the asymptotic
\begin{equation}\label{F-asymp}  F(x) = \exp\left( \left(\frac{C_0}{2}+o(1)\right) \Log_2^{1/2} x \Log_3 x \right)\end{equation}
as $x \to \infty$ (we allow the decay rates $o(1)$ here and in the sequel to depend on $C_0$), although the $o(1)$ error is somewhat oscillatory.

We then show

\begin{theorem}[Main theorem]\label{main-442}\ Let $C_0>0$ be a fixed constant.
  \begin{itemize}
    \item[(i)]
    There exists a set $A$ of natural numbers such that
  \begin{equation}\label{est-1} \sum_{n \in A: n \leq x} \frac{1}{n} \asymp F(x) \end{equation}
  for sufficiently large $x$, whose defect \eqref{defect} obeys the bound
  \begin{equation}\label{est-2}
    \E \gcd(\mathbf{n}, \mathbf{m}) - 1 \ll \exp(C_0^2)-1+o(1)
  \end{equation}
  as $x \to \infty$.
    \item[(ii)]  Conversely, if $x$ is sufficiently large, and $A$ is a set of natural numbers for which
\begin{equation}\label{est-2a}   \E \gcd(\mathbf{n}, \mathbf{m}) - 1 \leq C_0^2
\end{equation}
  holds for that choice of $x$, then
    \begin{equation}\label{est-3}  \sum_{n \in A: n \leq x} \frac{1}{n} \leq \exp( (C_0+o(1)) \Log_2^{1/2} x \Log_3 x ).
    \end{equation}
  \end{itemize}
\end{theorem}

It is easy to see that Theorem \ref{simplified-main} follows from Theorem \ref{main-442} combined with \eqref{F-asymp}.

In the small defect regime $C_0 \approx 0$, Theorem \ref{main-442} provides a reasonably tight bound: for small $0 < \eps \leq 1$, by setting $C_0 \coloneqq \log^{1/2}(1+\eps^2) = \eps+ O(\eps^3)$, Theorem \ref{main-442}(i) permits us to construct sets $A$ in which the defect \eqref{defect} is bounded by $\eps^2 + o(1)$ and
$$ \sum_{n \in A: n \leq x} \frac{1}{n} = \exp\left( \left(\frac{\eps}{2}+O(\eps^3)+o(1)\right) \Log_2^{1/2} x \Log_3 x \right)$$
while for any set $A$ with this defect bound, Theorem \ref{main-442}(ii) provides the upper bound
$$ \sum_{n \in A: n \leq x} \frac{1}{n} = \exp\left( \left(\eps+o(1)\right) \Log_2^{1/2} x \Log_3 x \right).$$
In particular, we see that the right-hand side of \cite[Proposition 2.1]{bergelson} can only exhibit a decay of $o(1)$ for sets $B$ with the growth rate $\sum_{n \in B: n \leq x} \frac{1}{n} \leq \exp( o(\Log_2^{1/2} x \Log_3 x))$, and that this upper bound cannot be improved.  The example already noted in \cite{bergelson} of the set of products of $k$ primes for fixed $k$ morally corresponds to the $C_0 = k / \Log_2^{1/2} x$ case of the theorem, although strictly speaking Theorem \ref{main-442}  does not apply in this regime, since we require $C_0$ to be independent of $x$.  In the large $C_0$ regime, there is an exponential gap in constants between (i) and (ii); we suspect that the bounds in (i) are closer to the truth, but we do not pursue\footnote{Will Sawin (private communication) has shown in fact that the bounds in (i) are essentially optimal; see Appendix \ref{sawin-app}.} this matter in this paper.

\subsection{Methods of proof}

For simplicity we restrict attention in this informal discussion to the case $C_0 \asymp 1$, avoiding some technical difficulties that emerge in the case of small $C_0$.

We first discuss Theorem \ref{main-442}(ii).  From the form \eqref{defect-alt} of the defect \eqref{defect}, we see from the hypotheses that for a typical prime $p$, the probability that $p$ divides $\mathbf{n}$ (or $\mathbf{m}$) has to be somewhat small; by linearity of expectation, this will imply that the expected number $\omega(\mathbf{n})$ of prime factors of $\mathbf{n}$ will be much smaller than the number $\Log_2 x$ predicted by the Hardy--Ramanujan (or Erd\H{o}s--Kac) law; in fact, a simple application of Cauchy--Schwarz and Mertens' theorem reveals that it is of size at most $(C_0+o(1)) \Log^{1/2}_2 x$.  However, standard arguments from the anatomy of integers show that the set of numbers with this few prime factors obey the bound \eqref{est-3}.  The rigorous proof of (ii) will largely follow these lines, though it also invokes Jensen's inequality at one point to handle the fact that $\omega(\mathbf{n})$ is only controlled in expectation, rather than uniformly.

For the proof of Theorem \ref{main-442}(i), it is natural in view of the proof of (ii) (as well as the aforementioned analysis of the set of numbers with $k$ prime factors for a fixed $k$) to look at the set of numbers $n$ with at most $C_0 \Log^{1/2}_2 n$ prime factors.  This essentially works, except for a number of technical complications arising from the fact that the quantity $C_0 \Log^{1/2}_2 n$ is usually not an integer, leading to the left-hand sides of \eqref{est-1}, \eqref{est-2} experiencing some unwanted ``discontinuous'' behavior whenever $C_0 \Log^{1/2}_2 x$ crosses an integer value, if the set $A$ is not chosen carefully. Nevertheless, the author was able to modify this construction to reduce the discontinuities to acceptable levels, though at the cost of introducing the somewhat unusual growth function \eqref{F-def}.  The case when $C_0$ is small presented particular challenges in this regard, being one of the primary reasons for the complexity of the construction.   A simplified\footnote{See also the first version of the paper, available at \url{https://arxiv.org/abs/2407.04226v1}, which handled the asymptotic limit $x \to \infty$ but only in the special case $C_0 = 1$, leading to an argument intermediate in complexity between the proof of Theorem \ref{main-442}(i) and the proof of Proposition \ref{simplified}.} version of the construction, which has a shorter proof but only applies for a single large range of $x$ rather than in the asymptotic limit $x \to \infty$, is given in Proposition \ref{simplified} below.

\section{Proof of Theorem \ref{main-442}}\label{442-sec}

We begin with Theorem \ref{main-442}(ii), which is easier to prove, and also helps motivate the construction for (i).  We may of course assume that the set $\{ n \in A: n \leq x\}$ is non-empty, as the claim is trivial otherwise.
By \eqref{defect-alt}, the hypothesis \eqref{est-2a} is equivalent to
\begin{equation}\label{sumd}
  \sum_{d>1} \phi(d) \P(d|\mathbf{n})^2 \leq C_0^2.
\end{equation}
Specializing $d$ to be prime, we conclude in particular that\footnote{All sums over $p$ are understood to be over primes.  This specialization to primes is inefficient in the large $C_0$ regime, and is the main cause for the exponential gap between upper and lower bounds in this case. In Appendix \ref{sawin-app} we give an argument of Sawin that avoids restriction to primes, and which removes this inefficiency.}
$$ \sum_{p \leq x} (p-1) \P( p | \mathbf{n})^2 \leq C_0^2.$$
By Cauchy--Schwarz and Mertens' theorem, we conclude that
\begin{equation}\label{cauchy}
   \sum_{p \leq x} \P( p | \mathbf{n}) \leq C_0 \left( \sum_{p \leq x}\frac{1}{p-1}\right)^{1/2} \leq (C_0+o(1)) \Log_2^{1/2} x.
\end{equation}
Thus, if $\omega(n) \coloneqq \sum_{p|n} 1$ denotes the number of prime factors of $n$, we see from linearity of expectation that
$$ \E \omega(\mathbf{n}) \leq (C_0+o(1)) \Log_2^{1/2} x.$$
By Jensen's inequality, this implies that
$$ \E \exp( - \omega(\mathbf{n}) \Log_3 x ) \geq \exp( - (C_0+o(1)) \Log_2^{1/2} x \Log_3 x ).$$
Unwinding the probabilistic notation and rearranging, we obtain
$$ \sum_{n \leq x; n \in A} \frac{1}{n} \leq \exp( (C_0+o(1)) \Log_2^{1/2} x \Log_3 x )
\sum_{n \leq x; n \in A} \frac{\exp( - \omega(n) \Log_3 x )}{n}.$$
If we crudely drop the constraint $n \in A$ on the right-hand side, one can bound that expression by the Euler product
$$\exp( (C_0+o(1)) \Log_2^{1/2} x \Log_3 x ) \prod_{p \leq x} \left(1 + O\left(\frac{\exp(-\Log_3 x)}{p}\right)\right)
$$
which we can bound in turn by
$$\exp\left( (C_0+o(1)) \Log_2^{1/2} x \Log_3 x +O\left(\sum_{p \leq x} \frac{\exp(-\Log_3 x)}{p}\right) \right).$$
By Mertens' theorem, the sum inside the $O()$ notation is bounded, thus we obtain \eqref{est-3} as required.

Now we turn to Theorem \ref{main-442}(i).  Given a set ${\mathcal P}$ of primes and a natural number $k$, let ${\mathcal P}^{[k]}$ denote the set of all products $p_1 \dots p_k$ of $k$ distinct primes in ${\mathcal P}$, and then write ${\mathcal P}^{[\leq k]} \coloneqq \bigcup_{0 \leq k' \leq k} {\mathcal P}^{[k']}$ to be the collection of products of at most $k$ distinct primes in ${\mathcal P}$.  We have the following basic calculation, that describes the logarithmic growth rate of ${\mathcal P}^{[\leq k]}$  when $k$ is not too large:

\begin{lemma}[Logarithmic growth rate]\label{lgr}  let ${\mathcal P}$ be a set of primes, and let $k \geq 0$ be an integer.  If $x \geq 1$ is such that
  \begin{equation}\label{kp}
    k \leq c \sum_{p \in {\mathcal P}: p \leq x} \frac{1}{p}
  \end{equation}
  for a sufficiently small absolute constant $c>0$, then
  \begin{equation}\label{log}
    \sum_{n \in {\mathcal P}^{[k]}: n \leq x} \frac{1}{n} \asymp \sum_{n \in {\mathcal P}^{[\leq k]}: n \leq x} \frac{1}{n} \asymp \frac{1}{k!} \left( \sum_{p \in {\mathcal P}: p \leq x} \frac{1}{p} \right)^k.
  \end{equation}
\end{lemma}

This lemma can be compared with \cite[Lemma 2.1]{gkm}, which in our notation essentially\footnote{Strictly speaking, the result of \cite[Lemma 2.1]{gkm} gives slightly sharper bounds, but also permits the primes $p_1, \dots, p_k \in {\mathcal P}$ to be repeated; however, it is not difficult to adapt the proof of that lemma to achieve \eqref{log-all}, basically by restricting attention to squarefree numbers.} asserts the bound
\begin{equation}\label{log-all} \sum_{n \in \bigcup_{k=0}^\infty {\mathcal P}^{[k]}: n \leq x} \frac{1}{n} \asymp \exp \left( \sum_{p \in {\mathcal P}: p \leq x} \frac{1}{p} \right)
\end{equation}
without the need for a hypothesis such as \eqref{kp}. Combining \eqref{log} with \eqref{log-all}, we thus see that if $\mathbf{n}$ is drawn randomly from $\{ n \in \bigcup_{k=0}^\infty {\mathcal P}^{[k]}: n \leq x\}$ using the logarithmic probability distribution \eqref{n-prob}, then the lower tail distribution of $\omega(\mathbf{n})$ behaves like the lower tail of a Poisson random variable with parameter $\sum_{p \in {\mathcal P}: p \leq x} \frac{1}{p}$, in the spirit of the Sath\'e--Selberg formula \cite{sathe}, \cite{selberg}.  Results similar to Lemma \ref{lgr} are available for natural density\footnote{We are indebted to the anonymous referee for this remark and references.}, and can in fact be extended to the regime where the small constant $c$ is replaced by $1+o(1)$; see \cite{halasz, halasz-2}, as well as \cite[Theorem 5.4.2]{mangerel} for a joint version involving multiple disjoint sets ${\mathcal P}_1, \dots, {\mathcal P}_m$.  On the other hand, for very large $k$ the asymptotics deviate from the naive Poisson model; see \cite{norton} and \cite{gkm} for some relevant calculations.

In our applications of Lemma \ref{lgr}, the hypothesis \eqref{kp} will be automatically satisfied, as the left-hand side will be comparable to $C_0 \Log^{1/2}_2 x$ and the right-hand side comparable to $\Log_2 x$, with $x$ assumed to be large.

\begin{proof}  The case $k=0$ is trivial, so we may assume $k \geq 1$, which by \eqref{kp} forces $x$ to be larger than any specified absolute constant.

The right-hand side of \eqref{log} increases at least geometrically in $k$ thanks to \eqref{kp}, so by the triangle inequality it suffices to establish the bound for ${\mathcal P}^{[k]}$, or equivalently that
\begin{equation}\label{log'}
  \sum_{n \in {\mathcal P}^{[k]}: n \leq x} \frac{k!}{n} \left( \sum_{p \in {\mathcal P}: p \leq x} \frac{1}{p} \right)^{-k} \asymp 1.
\end{equation}

We argue probabilistically.  Let ${\mathbf p}_1,\dots,{\mathbf p}_k$ be random elements of $\{ p \in {\mathcal P}: p \leq x\}$, chosen independently with distribution
$$ \P( {\mathbf p}_i = p_i ) \coloneqq \frac{1/p_i}{\sum_{p \in {\mathcal P}: p \leq x} \frac{1}{p}}.$$
Then each $n \in {\mathcal P}^{[k]}$ is attained by the random variable ${\mathbf p}_1 \dots {\mathbf p}_k$ with probability $\frac{k!}{n} \left( \sum_{p \in {\mathcal P}: p \leq x} \frac{1}{p} \right)^{-k}$.  Conversely, ${\mathbf p}_1 \dots {\mathbf p}_k$ will produce an element of ${\mathcal P}^{[k]}$ if and only if all the ${\mathbf p}_i$ are distinct.  Thus the bound \eqref{log'} is equivalent to
\begin{equation}\label{pkx}
   \P( {\mathbf p}_1 \dots {\mathbf p}_k \leq x \hbox{ and } {\mathbf p}_1,\dots, {\mathbf p}_k \hbox{ distinct}) \asymp 1.
\end{equation}
The upper bound is trivial, so we focus on the lower bound.
By the union bound, the probability that ${\mathbf p}_1,\dots,{\mathbf p}_k$ are not distinct is at most
$$ \sum_{1 \leq i < j \leq k} \P( {\mathbf p}_i = {\mathbf p}_j ).$$
But by construction we have
\begin{align*}
  \P( {\mathbf p}_i = {\mathbf p}_j ) &= \frac{\sum_{p \in {\mathcal P}: p \leq x} \frac{1}{p^2}}{\left(\sum_{p \in {\mathcal P}: p \leq x} \frac{1}{p}\right)^2} \\
  &\ll \frac{1}{(\sum_{p \in {\mathcal P}: p \leq x} \frac{1}{p})^2}
\end{align*}
so by \eqref{kp} the probability of the ${\mathbf p}_1,\dots,{\mathbf p}_k$ not being distinct is at most $O(c^2)$.  Next, the event ${\mathbf p}_1 \dots {\mathbf p}_k > x$ can be written as $\sum_{i=1}^k \log {\mathbf p}_i > \log x$, so by Markov's inequality and linearity of expectation the probability of this event is at most
$$ \frac{1}{\log x} \sum_{i=1}^k \E \log {\mathbf p}_i.$$
But by construction and Mertens' theorem we have
\begin{align*}
\E \log {\mathbf p}_i &= \frac{\sum_{p \in {\mathcal P}: p \leq x} \log p/p}{\sum_{p \in {\mathcal P}: p \leq x} 1/p} \\
  &\ll \frac{\log x}{\sum_{p \in {\mathcal P}: p \leq x} 1/p}
\end{align*}
so by \eqref{kp} the probability of the event ${\mathbf p}_1 \dots {\mathbf p}_k > x$ is at most $O(c)$.  Thus the left-hand side of \eqref{pkx} is $1 - O(c^2) - O(c)$, and the claim follows for $c$ small enough.
\end{proof}

As a quick application of this lemma, we can establish a ``baby'' version of Theorem \ref{main-442}(ii), in which we only control the behavior at a single large range  of $x$, rather than all large $x$:

\begin{proposition}[Baby version of theorem]\label{simplified} Let $C_0>0$ be a fixed constant, let $k$ be a natural number that is sufficiently large depending on $C_0$.  Then there exists a set $A$ of natural numbers (depending on $k$ and $C_0$) such that the bound
\begin{equation}\label{oop}
\sum_{n \in A: n \leq x} \frac{1}{n} = \exp\left( \left(\left(\frac{C_0}{2}+o(1)\right) \Log_2^{1/2} x \Log_3 x \right) \right)
\end{equation}
and the bound \eqref{est-2} both hold in the range
\begin{equation}\label{co}
 C_0 \Log_2^{1/2} x = (1+o(1)) k
 \end{equation}
where the $o(1)$ notation  is now with respect to the limit $k \to \infty$.
\end{proposition}

\begin{proof}
  We set $A \coloneqq {\mathcal P}^{[\leq k]}$, where ${\mathcal P}$ is the set of all primes.  From Lemma \ref{lgr} and Mertens' theorem we have
\begin{equation}\label{esta}
\sum_{n \in A; n \leq x} \frac{1}{n} \asymp \frac{1}{k!} \Log^k_2 x
\end{equation}
where
The claim \eqref{oop} then follows readily from Stirling's formula and \eqref{co}.

It remains to establish \eqref{est-2}.  By \eqref{esta} and \eqref{gcd-standard} (or \eqref{sum-square}), and bounding $\phi(d)$ by $d$, it suffices to show that
\begin{equation}\label{dsum-alt}
 \sum_{1 < d \leq x} d \left( \sum_{n \in A: n \leq x; d|n} \frac{1}{n} \right)^2 \ll (\exp(C_0^2)-1+o(1)) \frac{1}{(k!)^2} \Log^{2k}_2 x.
\end{equation}
Observe that in order for $d$ to contribute to the above sum, it must be the product of $\ell$ primes for some $1 \leq \ell \leq k$, in which case $n$ is of the form $n = dm$ for some $m \in {\mathcal P}^{[\leq k-\ell]}$.  We can thus upper bound the left-hand side of \eqref{dsum-alt} by
$$
\sum_{\ell=1}^k \sum_{d \in {\mathcal P}^{[\ell]}: d \leq x} \frac{1}{d}  \left( \sum_{m \in {\mathcal P}^{[\leq k-\ell]}: m \leq x/d} \frac{1}{m} \right)^2.$$
Applying Lemma \ref{lgr} and Mertens' theorem, we can upper bound this quantity by
$$
\ll \sum_{\ell=1}^k \frac{1}{((k-\ell)!)^2} \Log^{2k-2\ell}_2 x \sum_{d \in {\mathcal P}^{[\ell]}: d \leq x} \frac{1}{d}$$
and by a further application of Lemma \ref{lgr} and Mertens' theorem, we can bound this expression in turn by
$$
\ll \sum_{\ell=1}^k \frac{1}{\ell! ((k-\ell)!)^2} \Log^{2k-\ell}_2 x.$$
Note that
\begin{equation}\label{kal}
   \frac{1}{(k-\ell)!} = \frac{k^\ell}{k!} \prod_{j < \ell} (1-j/k) \leq \frac{k^\ell}{k!} 
\end{equation}
so we can bound the preceding expression by
$$
\ll \frac{1}{(k!)^2} \sum_{\ell=1}^k \frac{(k^2 / \Log_2 x)^{\ell}}{\ell!}.$$
Using the Taylor expansion of the exponential function, we can bound this by
$$
\ll \frac{1}{(k!)^2} \left( \exp(k^2 / \Log_2 x) - 1 \right).$$
which gives \eqref{est-2} thanks to \eqref{co}.
\end{proof}

Now we return to establishing the full strength of Theorem \ref{main-442}(ii), using a more complicated version of the above construction.  Define\footnote{This construction was arrived at after a lengthy trial and error process, centered around trying to construct a set $A$ where a typical element $n$ had about $C_0 \Log_2^{1/2} n$ prime factors.  In particular, the author experimented extensively with several possible probabilistic constructions for the set $A$ to deal with the technical issue of $C_0 \Log_2^{1/2} n$ not being an integer, but to the author's surprise, such constructions turned out to be inferior to the deterministic construction provided here.} $A$ to be the set of all numbers $n$ which are square-free, lie in the range $x_{k} < n \leq x_{k+1}$ for some $k \geq C_0$, and are of one of the following two forms:
\begin{itemize}
  \item[(i)] $n$ is the product of at most $k$ primes in the range $[1,x_{k}]$.
  \item[(ii)] $n$ is the product of at most $k$ primes in the range $[1,x_{k}^{\eps(n)}]$, times one prime in the range $(x_{k}, x_{k+1}]$, where the exponent $\eps(y) \in (0,1]$ for any (real) parameter $x_k < y \leq x_{k+1}$ is defined by the formula
\begin{equation}\label{log-def}
\eps(y) \coloneqq \left( \frac{h(y)}{h(x_{k+1})} \right)^{k/C_0^2}.
\end{equation}
\end{itemize}
The exponent $\eps()$ is related to the unusual form of the piecewise quadratic function \eqref{psi-def}, and is needed for technical reasons to smooth out a ``discontinuity'' at $x_k$ when $C_0$ is small; we recommend that the reader ignore this exponent for a first\footnote{In the case $C_0 \asymp 1$, one can replace $\eps(n)$ by $1$ and $\psi(x)$ in \eqref{psi-def} by $h(x)$, leading to some simplifications in the analysis below, with the case $C_0=1$ being particularly clean notationally; we refer the reader to the previous version of the paper at \url{https://arxiv.org/abs/2407.04226v1} for details.} reading.  The exponent $\epsilon(n)$ is quite small when $n$ is close to $x_k$, and increases steadily to $1$ when $n$ reaches $x_{k+1}$; thus the effect of this exponent is to damp out case (ii) in the former regime, but not in the latter.

We now establish \eqref{est-1}.  We may assume that $x_k < x \leq x_{k+1}$ for some large natural number $k$ (in particular, much larger than $C_0$); we then have
$$ \Log_2 x = \frac{k^2}{C_0^2} + h(x) - 1$$
with $1 \leq h(x) \ll k/C_0^2$.  We record the convenient asymptotic
\begin{equation}\label{jb}
  (\Log_2 x_k + O(1))^j, (\Log_2 x + O(1))^j =\left(\left(1+O\left(\frac{1}{k}\right)\right)k^2/C_0^2\right)^j \asymp k^{2j} C_0^{-2j}
\end{equation}
for any $j = O(k)$; this asymptotic will combine particularly well with Lemma \ref{lgr} and Mertens' theorem.
Note that
\begin{equation}\label{keps}
\begin{split}
\Log_2 x_k^{\eps(x)} &= \Log_2 x_k + \log \eps(x) \\
&= \frac{k^2}{C_0^2} \left(1 + \frac{1}{k} \log\left(\frac{h(x)}{h(x_{k+1})}\right) \right) \\
&= \frac{k^2}{C_0^2} \left(1 - O\left(\frac{\log k}{k}\right)\right)
\end{split}
\end{equation}
and hence
\begin{equation}\label{jb-2}
  \begin{split}
  (\Log_2 x_k^{\eps(x)} + O(1))^k &\asymp
  \left( 1+\frac{1}{k}\log\left(\frac{h(x)}{h(x_{k+1})}\right) \right)^k k^{2k} C_0^{-2k} \\
  &\asymp \frac{h(x)}{h(x_{k+1})} k^{2k} C_0^{-2k}.
  \end{split}
\end{equation}
The same argument also gives
\begin{equation}\label{jb-2'}
  (\Log_2 x_k^{\eps(x^{\exp(-h(x)/2)})} + O(1))^k \asymp \frac{h(x)}{h(x_{k+1})} k^{2k} C_0^{-2k},
\end{equation}
as well as
\begin{equation}\label{jb-3}
  \begin{split}
  (\Log_2 x_k^{\eps(x)} + O(1))^{k-\ell} &\asymp \left( \frac{h(x)}{h(x_{k+1})}\right)^{1-\ell/k} k^{2(k-\ell)} C_0^{-2(k-\ell)} \\
  &= \left(1 + O\left(\frac{\log k}{k}\right)\right)^\ell \frac{h(x)}{h(x_{k+1})} k^{2(k-\ell)} C_0^{-2(k-\ell)}\\
  &= (1 + o(1))^\ell \frac{h(x)}{h(x_{k+1})} k^{2(k-\ell)} C_0^{-2(k-\ell)}
  \end{split}
\end{equation}
for any $0 \leq \ell \leq k$.

Observe (from the non-decreasing nature of $\eps$ on $(x_k,x]$) that any $n \in A$ with $n \leq x$ either lies in ${\mathcal P}_{\leq x_k}^{[\leq k]}$, or else is the product of a prime in $(x_k, x]$ and an element of ${\mathcal P}_{\leq z}^{[\leq k]}$, where $z = z(x) \coloneqq x_k^{\eps(x)}$, and ${\mathcal P}_{\leq y}$ denotes the set of primes up to $y$ for any given real $y$.  Hence we may upper bound the left-hand side of \eqref{est-1} by
$$ \left( \sum_{n \in {\mathcal P}_{\leq x_k}^{[\leq k]}: n \leq x} \frac{1}{n} \right)
+ \left( \sum_{n \in {\mathcal P}_{\leq z}^{[\leq k]}: n \leq x} \frac{1}{n} \right) \sum_{x_k < p \leq x} \frac{1}{p}.$$
Applying Lemma \ref{lgr}, Mertens' theorem, and \eqref{jb}, \eqref{jb-2}, \eqref{h-def} this is bounded by
$$ \ll \frac{k^{2k} C_0^{-2k}}{k!} \left(1 + \frac{h(x)}{h(x_{k+1})}(h(x)+O(1))\right) \ll F(x)$$
as required, thanks to \eqref{F-def}.

Now we turn to the lower bound.  One contribution to this sum comes from those numbers of the form $pm$ for $x_{k-1/2} < p \leq x_k^{1/2}$ a prime and $m \in {\mathcal P}_{\leq x_{k-1}^{\eps(x_{k-1/2})}}^{[\leq k-1]}$ with $m \leq x_k^{1/2}$, where we extend the sequence $x_k$ to fractional $k$ in the obvious fashion.  Note that the product $pm$ uniquely determines both $p$ and $m$. By Lemma \ref{lgr}, Mertens' theorem, and \eqref{jb}, \eqref{jb-2}, and also noting that $h(x_{k-1/2}) \asymp h(x_k) \asymp k/C_0^2$, this contribution is
$$ \gg \frac{(k-1)^{2(k-1)} C_0^{-2(k-1)}}{(k-1)!} \sum_{x_{k-1/2} < p \leq x^{1/2}_k} \frac{1}{p}
\asymp \frac{(k-1)^{2(k-1)} C_0^{-2(k-1)}}{(k-1)!} \left(\frac{k}{C_0^2} + O(1)\right) \asymp \frac{k^{2k} C_0^{-2k}}{k!}$$
(cf. \eqref{F-cts}). By \eqref{F-def}, \eqref{F-flat}, this handles the case when $h(x) \leq k^{1/2}/C_0$, so we may now assume that $h(x) \geq k^{1/2}/C_0$, then by \eqref{F-def} we have
\begin{equation}\label{F-form}
  F(x) \asymp \frac{h(x)^2}{h(x_{k+1})} \frac{k^{2k} C_0^{-2k}}{k!}.
\end{equation}
In this case we consider a different contribution to the left-hand side of \eqref{est-1}, namely from those numbers of the form $pm$ for $x^{1/\exp(h(x)/2)} < p \leq \sqrt{x}$ and $m \in {\mathcal P}_{\leq z'}^{[\leq k]}$ with $m \leq \sqrt{x}$, where $z' \coloneqq x_k^{\eps(x^{1/\exp(h(x)/2)})}$.  By Lemma \ref{lgr}, Mertens' theorem, and \eqref{jb-2'}, this contribution is
$$ \gg \frac{h(x)}{h(x_{k+1})} \frac{k^{2k} C_0^{-2k}}{k!} \sum_{x^{1/\exp(h(x)/2)} < p \leq \sqrt{x}} \frac{1}{p} \asymp \frac{h(x)}{h(x_{k+1})} \frac{k^{2k} C_0^{-2k}}{k!} (h(x) - O(1)),$$
which is acceptable by \eqref{F-form} and the hypothesis $h(x) \geq k^{1/2}/C_0$.  This completes the proof of \eqref{est-1} in all cases.

  Now we establish \eqref{est-2}; this is a calculation similar to that used to esablish Proposition \ref{simplified}, but with more technicalities.  We again may assume that $x_k < x \leq x_{k+1}$ for some large natural number $k$.  By \eqref{est-1} and \eqref{gcd-standard} (or \eqref{sum-square}), and bounding $\phi(d)$ by $d$, it suffices to show that
\begin{equation}\label{dsum}
 \sum_{1 < d \leq x} d \left( \sum_{n \in A: n \leq x; d|n} \frac{1}{n} \right)^2 \ll (\exp(C_0^2)-1+o(1)) F(x)^2.
\end{equation}
Observe that in order for the inner sum to be non-zero, the natural number $1 < d \leq x$ must be one of the following two forms:
\begin{itemize}
  \item[(i)] $d$ is the product of $\ell$ primes in the range $[1,x_{k}]$ for some $1 \leq \ell \leq k$.
  \item[(ii)] $d$ is the product of $\ell$ primes in the range $[1,x_{k}]$ for some $0 \leq \ell \leq k$, times one prime $p$ in the range $(x_{k}, x_{k+1}]$.
\end{itemize}
First suppose we are in case (i).  Then the numbers $n$ that contribute to the inner sum  of \eqref{dsum} are of the form $dm$, where $m$ either lies in ${\mathcal P}_{\leq x_k}^{[\leq k-\ell]}$, or is an element of ${\mathcal P}_{\leq z}^{[\leq k-\ell]}$ times a prime in $(x_k, x]$, where we recall $z \coloneqq x_k^{\eps(x)}$.  We conclude the upper bound
\begin{align*}
\sum_{n \in A: n \leq x; d|n} \frac{1}{n} &\ll \frac{1}{d} \left( \sum_{m \in {\mathcal P}_{\leq x_k}^{[\leq k-\ell]}; m \leq x} \frac{1}{m} \right) \\
&\quad + \left( \sum_{m \in {\mathcal P}_{\leq z}^{[\leq k-\ell]}; m \leq x} \frac{1}{m} \right) \left(1 + \sum_{x_k \leq p < x} \frac{1}{p}\right)
\end{align*}
and hence by Lemma \ref{lgr}, Mertens' theorem, and \eqref{jb}, \eqref{jb-3}
$$ \sum_{n \in A: n \leq x; d|n} \frac{1}{n} \ll \frac{1}{d}  \frac{(1 + o(1))^\ell}{(k-\ell)!} k^{2k-2\ell} C_0^{-2k+2\ell}
\left(1 + \frac{h(x)}{h(x_{k+1})} (h(x)+O(1))\right).$$
From \eqref{psi-def} we have
$$1 + \frac{h(x)}{h(x_{k+1})} (h(x)+O(1)) \asymp  \psi(x).$$
The total contribution of case (i) to \eqref{dsum} is thus
\begin{equation}\label{best}
   \ll \psi(x)^2 \sum_{\ell=1}^k \frac{(1 + o(1))^\ell}{((k-\ell)!)^2} k^{4k-4\ell} C_0^{-4k+4\ell}
   \sum_{d \in {\mathcal P}_{\leq x_k}^{[\ell]}: d \leq x} \frac{1}{d}.
\end{equation}
We can estimate the inner sum using Lemma \ref{lgr}, Mertens' theorem, and \eqref{jb}, to obtain the bound
$$ \ll \psi(x)^2 \sum_{\ell=1}^k  \frac{(1 + o(1))^\ell}{\ell! ((k-\ell)!)^2} k^{4k-2\ell} C_0^{-4k+2\ell}.$$
Using \eqref{kal} as in the proof of Proposition \ref{simplified}, we can bound the preceding expression using \eqref{F-def} by
$$ \ll F(x)^2 \sum_{\ell=1}^\infty \frac{(C_0^2+o(1))^{\ell}}{\ell!}.$$
Performing the $\ell$ summation, we obtain the desired bound of $O((\exp(C_0^2)-1+o(1)) F(x)^2)$.

Now suppose we are in case (ii), which will be a lower order term (but only because of the ``flatness'' properties of $\psi$ near $x_k$).  Then the numbers $n$ that contribute to the inner sum of \eqref{dsum} are of the form $dm$ where $m$ lies in ${\mathcal P}_{\leq z}^{[\leq k-\ell]}$, so that
$$ \sum_{n \in A: n \leq x; d|n} \frac{1}{n} \ll \frac{1}{d} \left( \sum_{m \in {\mathcal P}_{\leq z}^{[\leq k-\ell]}; m \leq x} \frac{1}{m} \right)$$
and hence by Lemma \ref{lgr}, Mertens' theorem, and \eqref{jb-3}
$$ \sum_{n \in A: n \leq x; d|n} \frac{1}{n} \ll \frac{1}{d} \frac{h(x)}{h(x_{k+1})} \frac{ (1+o(1))^\ell}{(k-\ell)!} k^{2(k-\ell)}.$$
The total contribution of case (ii) to \eqref{dsum} is thus
$$ \ll \left( \frac{h(x)}{h(x_{k+1})}\right)^2 \sum_{\ell=0}^k \frac{(1+o(1))^\ell}{((k-\ell)!)^2}  k^{4(k-\ell)}\left(\sum_{d \in {\mathcal P}_{\leq x_k}^{[\ell]}: d \leq x} \frac{1}{d}\right) \sum_{x_k \leq p < x} \frac{1}{p}.$$
By Mertens' theorem, and Lemma \ref{lgr}, Mertens' theorem, and \eqref{jb}, this is bounded by
$$ \ll (h(x)+O(1)) \left( \frac{h(x)}{h(x_{k+1})}\right)^2 \sum_{\ell=0}^k\frac{ (1+o(1))^\ell }{((k-\ell)!)^2 (\ell!)^2} k^{4k-2\ell}.$$
From \eqref{psi-def} we have
\begin{align*}
   (h(x)+O(1)) \left( \frac{h(x)}{h(x_{k+1})}\right)^2 &\ll h(x_{k+1})^{-1/2}
\left( 1 + \frac{h(x)^2}{h(x_{k+1})}\right)^{3/2} \\
&\ll h(x_{k+1})^{-1/2}
\left( 1 + \frac{h(x)^2}{h(x_{k+1})}\right)^2\\
& = h(x_{k+1})^{-1/2} \psi(x)^2 \\
&= o(\psi(x)^2);
\end{align*}
from this and \eqref{kal} we may bound the preceding expression by
$$ \ll o\left( \psi(x)^2 \frac{k^{4k}}{(k!)^2} \sum_{\ell=0}^\infty \frac{(C_0^2+o(1))^{\ell}}{\ell!}\right).$$
Performing the $\ell$ sum and using \eqref{F-def}, we obtain a bound of $o(F(x)^2)$.  Combining this with the case (i) analysis,  and the claim \eqref{est-2} follows.

\vskip20pt\noindent {\bf Acknowledgements.}

The author is supported by NSF grant DMS-1764034.  The author is indebted to Joel Moreira for pointing out the connection to the preprint \cite{bergelson}, Fabrice Imparato and Will Sawin for corrections, Marc Munsch for pointing out the reference \cite{dmt}, and the anonymous referee for helpful suggestions and corrections.

\appendix

\section{An argument of Sawin}\label{sawin-app}

In this appendix we give an argument of Will Sawin (private communication), which was provided to us after the initial release of this paper, and which gives an improved version of Theorem \ref{main-442}(ii) that matches the asymptotics in Theorem \ref{main-442}(i) to leading order:

\begin{theorem}[Improved upper bound]  If $C_0>0$ is a fixed constant, $x$ is sufficiently large, and $A$ is a set of natural numbers such that
\begin{equation}\label{est-2a-w}   \E \gcd(\mathbf{n}, \mathbf{m}) - 1\leq e^{C_0^2} - 1
\end{equation}
holds for that choice of $x$, then
    \begin{equation}\label{est-3-w}  \sum_{n \in A: n \leq x} \frac{1}{n} \leq \exp\left( \left(\frac{C_0}{2}+o(1)\right) \Log_2^{1/2} x \Log_3 x \right).
    \end{equation}
  \end{theorem}

\begin{proof} Similar to \eqref{sumd}, the hypothesis \eqref{est-2a-w} can be written as
$$ \sum_{d \geq 1} \phi(d) \P(d|\mathbf{n})^2 \leq e^{C_0^2}.$$
In contrast to the proof of Theorem \ref{main-442}(ii), we do not restrict $d$ to be prime.  Instead, we use Cauchy--Schwarz to observe for any $\delta>0$ that 
\begin{align*}
   \E (1 + \delta)^{\omega(\mathbf{n})} &= \sum_{d \leq x} \delta^{\omega(d)} \mu^2(d) \P(d|\mathbf{n}) \\
   &\leq e^{C_0^2/2} \left( \sum_{d \leq x} \frac{\delta^{2\omega(d)} \mu^2(d) }{\phi(d)} \right)^{1/2}.
\end{align*}
We can bound
\begin{align*}
 \sum_{d \leq x} \frac{\delta^{2\omega(d)} \mu^2(d) }{\phi(d)} 
&= \sum_{k=0}^\infty \delta^{2k} \sum_{d \in {\mathcal P}_{\leq x}^{[k]}} \frac{1}{\phi(d)} \\
&\leq \sum_{k=0}^\infty \frac{\delta^{2k}}{k!} \left( \sum_{p \leq x} \frac{1}{p-1} \right)^k \\
&= \exp\left( \delta^2 \sum_{p \leq x} \frac{1}{p-1} \right) \\
&= \exp\left( (1+o(1)) \delta^2 \Log_2 x \right).
\end{align*}
thanks to Mertens' theorem, where ${\mathcal P}_{\leq x}$ denotes the set of primes up to $x$.  On setting $\delta \coloneqq C_0 / \Log_2^{1/2} x$ we obtain
$$\E (1 + C_0 / \Log_2^{1/2} x)^{\omega(\mathbf{n})} \leq e^{(1+o(1))C_0^2}.$$
In particular, from Markov's inequality, we see for any fixed $\eps>0$ that
$$ \P\left((1 + C_0 / \Log_2^{1/2} x)^{\omega(\mathbf{n})} \leq e^{(1+\eps + o(1))C_0^2} \right) \gg_\eps 1.$$
The event
$$ (1 + C_0 / \Log_2^{1/2} x)^{\omega(\mathbf{n})} \leq e^{(1+\eps + o(1))C_0^2}$$
implies that
$$ \omega(\mathbf{n}) \leq (C_0+C_0 \eps + o(1)) \Log_2^{1/2} x$$
and hence by the definition \eqref{n-prob} of $\mathbf{n}$ one has
$$  \sum_{n \in A: n \leq x} \frac{1}{n}  \ll_\eps \sum_{n \leq x: \omega(n) \leq (C_0+C_0 \eps + o(1)) \Log_2^{1/2} x} \frac{1}{n}.$$
By Lemma \ref{lgr} and Mertens' theorem, one has
$$
\sum_{n \leq x: \omega(n) \leq (C_0+C_0 \eps + o(1)) \Log_2^{1/2} x} \frac{1}{n} \asymp \frac{1}{k!} (\Log_2 x + O(1))^k 
$$
for some $k = (C_0+C_0 \eps + o(1)) \Log_2^{1/2} x$.  Using Stirling's formula, we conclude that
$$ \sum_{n \in A: n \leq x} \frac{1}{n}
\ll_\eps \exp\left( \left(\frac{C_0}{2} + \frac{C_0 \eps}{2} + o(1)\right) \Log_2^{1/2} x \Log_3 x \right);$$
sending $\eps$ slowly to zero, we obtain the claim.
\end{proof}

\end{document}